\documentclass[letterpaper, 10 pt, conference]{ieeeconf}  
\usepackage[utf8]{inputenc}
\usepackage{amsmath}
\usepackage{amsfonts}
\usepackage{mathtools}
\usepackage{float}
\usepackage{graphicx}

\usepackage{mathtools}  
\usepackage{algorithm}
\usepackage{algorithmic}
\usepackage{xcolor}
\mathtoolsset{showonlyrefs=false}

\IEEEoverridecommandlockouts                              

\makeatletter
\def\endthebibliography{%
  \def\@noitemerr{\@latex@warning{Empty `thebibliography' environment}}%
  \endlist
}
\makeatother

\title{\LARGE \bf
Platoon Coordination and Leader Selection in Mixed Transportation Systems via Dynamic Programming
}

\author{Ying Wang$^{1}$, Ting Bai$^{2}$, and Andreas A. Malikopoulos$^{2}$,~\IEEEmembership{Senior Member, IEEE}
\thanks{*This work was supported by Wallenberg AI, Autonomous Systems and Software Program (WASP), funded by Knut and Alice Wallenberg Foundation. This research was also supported in part by NSF under Grants CNS-2401007, CMMI-2348381, IIS-2415478, and in part by MathWorks.}
\thanks{Y. Wang is with the Division of Decision and Control Systems,
	KTH Royal Institute of Technology, 10044 Stockholm, Sweden. E-mail: {\tt\small yinwang@kth.edu}}
\thanks{T. Bai, and A. A. Malikopoulos are with the School of Civil $\&$ Environmental Engineering, Cornell University, Ithaca, New York, USA. E-mails: \{{\tt\small tingbai, amaliko\}@cornell.edu}}
}

\begin{document}
\maketitle
\thispagestyle{empty}
\pagestyle{empty}

\begin{abstract}
With the growing penetration of electric trucks, freight transportation is transitioning toward a mixed system comprising both fuel-powered and electric trucks. Enhancing the formation of truck platoons in such a heterogeneous environment presents new challenges. This paper investigates the hub-based platoon coordination problem in a mixed-truck fleet, where the focus is to optimize trucks' waiting times, charging amounts for electric trucks, and platoon leader assignments. The objective is to maximize the overall platoon revenue of the fleet while accounting for the associated waiting and charging costs. We formulate the problem as a centralized optimization problem and present a dynamic programming approach to compute its suboptimal solution efficiently. The proposed method operates in polynomial time, ensuring scalable computational efficiency. Simulation studies involving 1,000 trucks traveling between two hubs in Sweden demonstrate the effectiveness and scalability of the approach. 
\end{abstract}

\vspace{-2pt}\section{Introduction}
Truck platooning has emerged as a promising technology in the context of emerging mobility and freight transportation systems. By utilizing vehicle-to-vehicle communication and automated control strategies, platooning enables heavy-duty trucks to travel in tightly spaced formations, thereby reducing aerodynamic drag. This reduction leads to significant improvements in fuel efficiency and associated environmental benefits~\cite{Bai2022-pz}. 
To be able to fully enjoy platooning benefits, trucks with diverse pre-planned routes and schedules need to be coordinated to form platoons safely and seamlessly. Existing coordination strategies facilitate platoon formation through a variety of means, including route planning to maximize shared road segments~\cite{larson2013coordinated}, 
adjusting vehicle speeds to allow for merging into platoons on route~\cite{guo2018fuel}, and scheduling departure or waiting times at logistics hubs~\cite{johansson2022platoon,9683080}, with drivers' hours-of-service regulations incorporated~\cite{9993403}. Although these works contribute significantly to improving the platoon formation rate in road networks, the models and solutions developed are suitable only for fuel-powered trucks, limiting their applicability in mixed transportation systems.

The electrification of road freight transportation has emerged as a global trend aimed at reducing carbon emissions, addressing energy shortages, and advancing the transition to a greener transportation system~\cite{jahangir2021road,bakker2025strategic}. As electric trucks (ETs) gain increasing market penetration, conventional freight networks are evolving into mixed transportation systems where fuel-powered trucks (FTs) and ETs coexist. Unlike FTs, ETs have limited battery capacity and often require mid-route charging schedules~\cite{10147895} to complete long-haul deliveries on time, especially given limited charging resources~\cite{10932686}. This introduces new scheduling challenges and potential disruptions to coordinated platooning. In addition, energy-saving benefits from platooning differ between FTs and ETs, making leader selection a critical decision.

While there have been growing efforts in platoon coordination for mixed traffic scenarios, efficient strategies that jointly address charging planning, leader selection, and waiting time coordination within a unified framework remain scarce. The authors in \cite{lee2021optimal} presented an eco-friendly strategy to optimize platoon configurations for a heterogeneous fleet of electric vehicles, with the target of maximizing the platooning advantages. However, their approach focuses on optimal platoon formation without taking into account the vehicles' charging requirements. More recently, \cite{ALAM2023104009} introduced a mixed-integer linear programming model to co-optimize charging strategies and platooning for long-haul electric freight vehicles. Besides, \cite{scholl2023platooning} developed an E-platooning scheme that coordinates vehicle charging with the waiting times needed to form platoons with potential partners. While these studies successfully address platoon formation among ETs, they overlook opportunities arising in mixed fleets with both ETs and FTs. The authors in~\cite{yao2020managing} proposed a decentralized method to assign leader and follower roles in mixed-traffic platoons. Nevertheless, their framework assumes that all trucks form one platoon and depart simultaneously from the station, with optimization limited to only the leader's role. 

This paper aims to develop a hub-based platoon coordination strategy in a mixed transportation system composed of both FTs and ETs. We focus on a fleet of trucks that arrive at a hub with a random arrival process. Based on the battery levels of the ETs upon their arrival, our coordination scheme enables charging time scheduling for ETs, waiting time scheduling for ETs and FTs, and leader selection for each platoon, where the target is to maximize the overall platooning profit of the fleet. The main contributions are: (i) We investigate the platoon coordination problem
in a mixed transportation system, taking into account different models of FTs and ETs, and their scheduling constraints. (ii) We propose a dynamic programming approach to address the platoon coordination problem approximately and efficiently. The solution ensures that ETs operate safely while maximizing the overall platoon utility in a discretized solution space. 


\vspace{-3.5pt}\section{System Models}
\label{sec:Problem Setting}

We consider a road network consisting of hubs and road segments connecting hubs, where hubs, in practice, represent parking areas or service centers. We consider $N\!\in\!\mathbb{N}$ trucks at a hub $A$ having the road segment overlapping on routes aim to travel in platoons heading to the same connected hub $B$. The distance between the two hubs is denoted as $d\!\in\!\mathbb{R}_+$. For modeling simplicity, we assume that all trucks travel at an identical and constant speed. Let $\mathcal{N}\!:=\!\{1,2,\dots,N\}$ be the set of all trucks, consisting of FTs and ETs, whose sets are denoted as $\mathcal{F}$ and $\mathcal{E}$, respectively, i.e., $\mathcal{N}\!=\!\mathcal{F}\!\cup\!{\mathcal{E}}$. 

Trucks arrive at hub $A$ at different times, and we use $\tau_i^0$ to denote the arrival time of truck $i\!\in\!\mathcal{N}$. To form platoons with other trucks, each truck needs to optimize its waiting time at the hub based on its arrival time and the planned departure times of others. For ETs, in addition to waiting times, the charging periods need optimization to be able to reach the next hub $B$ with sufficient battery guarantees. To model the platoon coordination problem, we first introduce the dynamic models for each type of truck. 

\vspace{-4.5pt}\subsection{Departure Times}\vspace{-3pt}
Forming platoons requires all participating trucks to depart from hub $A$ at the same time. To describe this constraint, let $\tau_i^w\!\in\!\mathbb{R}_+$ denote the waiting time of truck $i\!\in\!\mathcal{N}$ and $\tau_i^c\!\in\!\mathbb{R}_+$ denote its charging time if $i\!\in\!\mathcal{E}$. Note that for an ET, the waiting time $\tau_i^w$ excludes the charging time and represents only the time spent waiting to join a platoon. Thus, for any truck $i\!\in\!\mathcal{N}$, its departure time from hub $A$ is determined as
\begin{align}
  \tau_i^{d}= 
\begin{cases}
   \tau_i^0 + \tau_i^w, \quad &{\rm if}  \ \ i \!\in\!\mathcal{F}, \\
   \tau_i^0  + \tau_i^c + \tau_i^w,  \quad &{\rm if}  \ \ i \in \mathcal{E}.
\end{cases}\label{Eq.1}
\end{align}

Due to limited space at hub $A$, we assume that each truck's stay is bounded by a planning horizon, denoted as $T\!\in\!\mathbb{N}_+$. Specifically, $\tau_i^d\!\leq\!{T}$ for $i\!\in\!\mathcal{N}$, where the start of the planning is set as $0$. Furthermore, we assume that the charging resources at hub $A$ are sufficient, ensuring that no extra waiting time arises due to charging congestion.   


\vspace{-4.5pt}\subsection{Charging and Discharging Models}\vspace{-3pt}
\label{sec:SoC_dym}
The battery level of each ET, described by the state-of-charge (SoC) given in percentage, is critical to ensure the safe and reliable operation of the truck throughout its journey. The dynamics of the SoC are determined by the charging and discharging models of the truck, as we introduce below.
\subsubsection{Charging model}
Throughout this paper, we assume that all ETs in the fleet have the same charging model. Let $s_i^0$ denote the SoC when truck $i$ arrives at hub $A$ and $s_i^d$ be the SoC when leaving the hub. The charging process of each ET is approximately modeled by a linear process 
\begin{equation}
    s_i^d = s_i^0 + r_i  \tau_i^c,
    \label{eq:charge_eq}
\end{equation}
where \( r_i\!\in\!\mathbb{R}_+ \) denotes the charging rate of truck $i$ (measured in the increase in SoC percentage per unit time). As defined previously, $\tau_i^c$ denotes the charging period. Note that, $s_i^d$ is constrained by the truck's battery capacity and must be sufficient to cover the energy required to reach the destination hub $B$ with safety guarantees during the journey. 

\subsubsection{Discharging model} For any ET $i\!\in\!\mathcal{E}$, let $v_i\!\in\!\mathbb{R}_+$ be its basic discharging rate (measured in the decrease in SoC percentage per unit travel distance), which represents the rate of energy consumption when the truck travels alone. Recall that $d$ denotes the distance between hub $A$ and hub $B$. Thus, the SoC of truck $i$ arriving at hub $B$, denoted as $s_i^a$, is
\begin{equation}
     s_i^a=s_i^d - \beta_i v_i  d, 
    \label{eq:discharge_eq}
\end{equation}
where $\beta_i$ denotes a coefficient that varies depending on the role of the ET when it is part of a platoon, i.e., acting as a leader or a follower. Since trucks traveling as followers in a platoon experience lower energy consumption, resulting in a reduced discharging rate, while those traveling as leaders receive significantly fewer energy savings. We determine the value of $\beta_i$, $i\!\in\!\mathcal{E}$ according to the following rule:
  \begin{itemize}
      \item If truck $i$ travels \emph{alone} or as a platoon \emph{leader}, $\beta_i\!=\!1$.
      \item If truck $i$ travels as a platoon \emph{follower}, $\beta_i\!=\!0.82$~\cite{tsugawa2016review}.
  \end{itemize}
  
To ensure safe operation, the SoC of each truck is no less than $s_i^{\rm safe}$, which denotes the minimum battery level allowed. As it requires that $s_i^a\!\geq\!{s_i^{\rm safe}}$, $s_i^d$ given above is confined by
\begin{equation} \label{eq:tau_c_constraint}
    s_i^{\rm safe}+\beta_i v_i d\leq s_i^d
    \leq s^{\max}_i, 
\end{equation}
where $s_i^{\rm max}$ represents the battery capacity of truck $i$, and $\beta_i$ is determined by the role of the truck in the platoon.

\vspace{-3.5pt}\section{Platoon Coordination Problem}\label{Section III}\vspace{-3.5pt}

Hub-based platoon coordination synchronizes the departure of trucks arriving at a shared hub at different times, enabling them to travel in platoons to a common destination. 
We consider a scenario where all trucks, including both FTs and ETs, belong to a single fleet.
The goal is to maximize fleet utility by optimally coordinating \textit{platoon departure times, ET charging schedules, and leader role assignments within platoons}.
The utility is defined as platoon profit minus its charging and waiting costs, as introduced below.

\subsubsection{Loss}
\label{sec:loss}

The loss associated with platoon formation consists of waiting loss and charging loss. Note that refueling costs for FTs and charging costs for ETs are excluded from these losses, as they are inherent travel costs that remain constant regardless of platoon formation.

The \textit{waiting loss} arises from additional labor expenses due to drivers' idle times for the extra waiting. This cost per unit of waiting time for each truck is denoted by \(\epsilon^w\!\in\!\mathbb{R}_+\). The \textit{charging loss}, represented by \(\epsilon^c\!\in\!\mathbb{R}_+\), is set lower than the waiting loss per unit time, reflecting a flexibility benefit since ETs can charge instead of waiting passively, thus increasing their SoC and operational flexibility upon arrival at the next hub B.
Formally, the total loss of the fleet \(\mathcal{N}\) is given by
\begin{align}      
    & L(\mathcal{N}) = \sum_{i \in \mathcal{E}}  (\epsilon^c \tau^c_i\!+\!\epsilon^w \tau^w_i) + \sum_{i \in \mathcal{F}} \epsilon^w \tau^w_i.
    \label{eq:loss_per_p}
\end{align}

\subsubsection{Profit}
\label{sec:profit}

Platoon profit comes from reduced fuel or electricity consumption by follower trucks. Field experiments show negligible savings for platoon leaders; thus, their contributions are ignored.
Profit over a distance \(d\) for each FT or ET follower is denoted by \(\xi^F\!\in\!\mathbb{R}_+\) or \(\xi^E\!\in\!\mathbb{R}_+\), respectively. 
Thus, the overall profit, obtained by summing the platoon profits at each time $t$ (with $\mathcal{P}_t$ denoting the platoon departing at time $t$, possibly empty if no vehicles depart then), is 
\begin{equation}
\label{eq:profit}
R(\mathcal{N}) = \sum_{t=1}^T R(\mathcal{P}_t)= \sum_{t=1}^T(\xi^E n^{E,f}_t+\xi^F n^{F,f}_t),
\end{equation}
where  $n^{E,f}_t$, $n^{F,f}_t$
denote the numbers of ET, FT followers in $\mathcal{P}_{t}$, respectively, decided by the selection of the leader type.

\subsubsection{Utility}
\label{sec:utility2}
By designing $\tau_i^c$, $\tau_i^w$ and leader roles, the fleet utility $J(\mathcal{N})$, obtained by $
J(\mathcal{N}) =R(\mathcal{N})-L(\mathcal{N})
$, is to be optimized.
Unfortunately, the search solution space grows exponentially in the fleet size $N$ and horizon length $T$, rendering a computationally prohibitive exhaustive exploration. This motivates the discretization approach and the structured platoon formation scheme presented in the following section.

\vspace{-3pt}\section{Dynamic Programming Solution}
\label{sec:Dynamic Programming}\vspace{-2pt}

In this section, we begin with a discretization scheme by introducing each vehicle’s earliest departure time, thereby defining a solution space for the dynamic programming (DP) search. We then present the DP utility function and its recurrence relation, and finally show that this DP‐based strategy produces a sub‐optimal solution in polynomial time.

\vspace{-4.5pt}\subsection{Discretize Solution Space}\vspace{-3pt}
\label{sec:ordering}

The primary computational challenge arises from the extensive solution space.
A possible simplification, in scenarios that involve only FTs, is to discretize the solution space by restricting the set of feasible platoon departure times $\mathcal{T}\!\coloneqq\!\{1, 2,\dots, T\}$ to the default arrival times $\{\tau_i^0:  i \in \mathcal{N}\}$ 
without compromising optimality. This simplification is justified by the observation that selecting any departure time between two consecutive default arrival times $\tau_i^0$ and $\tau_{i+1}^0$ would only increase the waiting loss without improving the platoon profit.
However, in mixed transportation systems, the main obstacle to directly applying this discretization to ETs is that they may have insufficient SoCs to reach the destination, necessitating additional charging time for safe travel.

Consider that each ET must be sufficiently charged to safely reach the next hub, at least in the most energy efficient scenario, i.e., as a platoon follower. We precompute the \textit{earliest departure time} for each ET based on the lower bound of the departure SoC specified in \eqref{eq:tau_c_constraint}. Given the initial SoC \( s^0_i \), road length \( d \), and energy consumption rate \( \beta_i v_i\), the minimum required SoC for safe departure is \(\beta_i v_i d + s_i^{\rm safe} \) with the follower energy coefficient $\beta_i\!=\!0.82$. Thus, we determine the minimum charging time to ensure safety as
\begin{equation}
    \label{eq:min_charge_time}
    \tau_i^{c, \min} = \{0, \ ( \beta_i v_i  d + s_i^{\rm safe} - s^0_i)/r_i \},
\end{equation}
with which the minimum departure SoC of the ET, denoted by $s_i^{d, \min}$, can be obtained with the linear charging model in \eqref{eq:charge_eq}.
Thus, the earliest departure time, accounting for the necessary minimal charging time of an ET is defined as
\begin{equation}
\label{eq:earliest departure}
   \tau_i^{\Delta}:= 
\begin{cases}
    \tau_i^0 + \tau_i^{c, \min}, \quad &{\rm if}  \ \ i \in \mathcal{E},\\
    \tau_i^0, \quad &{\rm if}  \ \ i \in \mathcal{F},
\end{cases}
\end{equation}
where \( \tau_i^0 \) denotes the default arrival time at hub A. For FTs, \( \tau_i^{\Delta}\) remains equal to their original arrival time.

Finally, all trucks are sorted by these earliest departure times, creating a discretized solution space $\{\tau_i^{\Delta}:i\in \mathcal{N}  \}$ for the subsequent DP search.
The complete discretization procedure is summarized in Algorithm~\ref{alg:1}.

\begin{algorithm}[t!]
\caption{Discretize Solution Space}
\label{alg:1}
\textbf{Input:}  Road length \(d\); for \(i \in \mathcal{N}\), arrival time \(\tau_i^0\), follower energy coefficient \(\beta_i\); for \(i \in \mathcal{E}\), initial SoC \(s_i^0\), charging rate \(r_i\), discharging rate \(v_i\), safe operating SoC \(s_i^{\rm safe}\).\\[1mm]
\textbf{Output:} The ordered set of earliest departure times $\{\tau_i^{\Delta}: i \in \mathcal{N}\}$ and, for ETs, the associated minimum departure SoC set
 $\{s_i^{d, \min}: i \in \mathcal{E}\}$
\begin{algorithmic}[1]
    \FOR{each vehicle \(i \in \mathcal{N}\)}
        \IF{\(i \in \mathcal{E}\)}
            \STATE Obtain the minimal charging time $\tau_i^{c, \min}$ with \eqref{eq:min_charge_time};
            \STATE Set the earliest departure time and departure SoC:
            $
                \tau_i^{\Delta} \gets \tau_i^0+\tau_i^{c, \min}, \quad s_i^{d, \min}  \gets s_i^0 + r_i \tau_i^{c, \min}.
            $
        \ELSE
            \STATE Set the earliest departure time for FTs:
            $\tau_i^{\Delta} \gets \tau_i^0.$ 
        \ENDIF
    \ENDFOR
    \STATE Sort all vehicles by their \(\tau_i^{\Delta}\) values in ascending order.
\end{algorithmic}
\end{algorithm}

\vspace{-4.5pt}\subsection{Platoon Utility Function}
\label{sec:dp cost}\vspace{-3pt}

Given the discretized solution space, the set of potential platoons~$\{\mathcal{P}_t: t\in \mathcal{T}\}$~is restricted to~$\{\mathcal{P}_{\tau_i^{\Delta}}: i \in \mathcal{N}\}$.
That is the set of feasible platoon departure times is reduced from the full range $\mathcal{T}=\{1,2,\dots,T\}$ to its subset \(\{\tau_i^{\Delta}:i \in \mathcal{N} \}\).

In our framework, we adopt a platoon formation scheme that \textit{restricts platooning to adjacent trucks}, defined as those that are consecutive in the order of their earliest departure times. This design is grounded in the observation that each follower truck receives a fixed profit regardless of the specific platoon it joins. Thus, there is no incentive for trucks to delay their departure in order to form platoons with non-adjacent vehicles. Instead, such delays would unnecessarily increase total waiting costs without providing any additional economic benefit. By limiting platoons to adjacent trucks, the scheme ensures efficient coordination while maintaining solutions that remain close to optimal overall platoon utility.



Formally, we define a platoon as a group of consecutive trucks: for truck \( i \in \mathcal{N} \), a platoon of size \( n \) comprising $i$ and its \( n-1 \) predecessors is represented as $\mathcal{P}_{i}^{n,m} \!=\! \big\{i-n+1, i-n+2, \ldots, i\big\}
$, where the other superscript \(m \!\in\! \{\texttt{E}, \texttt{F}\}\) indicates that the platoon leader is an ET if \(m\!=\!\texttt{E}\), or an FT if \(m\!=\!\texttt{F}\).
Note that this definition implies that the platoon is scheduled to depart at the earliest departure time \(\tau_i^{\Delta}\) of the truck \( i \).
The platoon utility, denoted as \(J(\mathcal{P}_i^{n,m})\), is given by
\begin{equation}
J(\mathcal{P}_i^{n,m}) =R(\mathcal{P}_i^{n,m})-L(\mathcal{P}_i^{n,m}).
\label{eq:specific_J}    
\end{equation}
The profit $R(\mathcal{P}_i^{n,m})$, defined in \eqref{eq:profit}, is calculated as
\begin{equation}
\label{eq:new profit}
R(\mathcal{P}_i^{n,m}) =
\begin{cases}
    \xi^F(|\mathcal{F}_i|\!-\!1) + \xi^E  |\mathcal{E}_i|    , & \text{if } m= \texttt{F},\\
   \xi^F |\mathcal{F}_i|  +  \xi^E (|\mathcal{E}_i|\!-\!1), & \text{if } m= \texttt{E},
\end{cases} 
\end{equation}
where $\mathcal{E}_i$ and $\mathcal{F}_i$ represent the sets of ETs and FTs in the platoon $\mathcal{P}_i^{n,m}$.
Their cardinalities correspond to the number of trucks in each set.
The loss $L(\mathcal{P}_i^{n,m})$ in \eqref{eq:specific_J}, accounting for both the waiting and charging cost, is expressed as
\begin{equation}
    L(\mathcal{P}_i^{n,m}) =   
      \sum_{j\in \mathcal{E}_i } \big(  \epsilon^c \tau_j^c+ \epsilon^w  \tau_j^w \big) 
     +\sum_{ j\in \mathcal{F}_i } \epsilon^w \tau_j^w, \label{eq:dp_loss}
\end{equation}
where for each FT $j\in \mathcal{F}_i$, $\tau_j^w=\tau_i^{\Delta} - \tau_j^{\Delta}$ since the departure time of $\mathcal{P}_i^{n,m}$ is $\tau_i^{\Delta}$. Similarly, for each ET $j\in \mathcal{E}_i$, $\tau_j^c + \tau_j^w=\tau_i^{\Delta} - \tau_j^{\Delta}$, where both the waiting time $\tau_j^w$ and the charging duration $\tau_j^c$ remain to be determined. 
Fortunately, as explained in Section \ref{Section III}, the monetary cost per unit of waiting time, $\epsilon^{w}$, is higher than the cost per unit of charging time, $\epsilon^{c}$.
Thus, the optimal solution naturally prioritizes charging over waiting, leading to a fixed charging policy.
In other words, ETs charge until they either reach the maximum capacity, $s_j^{\max}$, or the platoon's departure time, $\tau_i^{\Delta}$, whichever comes first, and any remaining time is then considered as waiting time. Thus, for ET \(j\!\in\! \mathcal{E}_i\), the total charging time, including the necessary minimal charging time $\tau_j^{c, \min}$ defined in \eqref{eq:min_charge_time}, is obtained as
\begin{equation}
    \tau_j^c =\tau_j^{c, \min} +   \{ (s_j^{\max}\!-\!s_j^{d,\min})/{r_j}, \  \tau_i^{\Delta}\!-\!\tau_j^{\Delta} \}. \label{eq:charge_for_ev}  
\end{equation}
After establishing the fixed charging strategy, both the charging and waiting time of each ET $j$ in the platoon $\mathcal{P}_i^{n,m}$ are determined, contributing to the following platoon loss:
\begin{equation}
    \begin{aligned}
    L(\mathcal{P}_i^{n,m})  
    \!=\!\sum_{j\in \mathcal{E}_i } \!\big(  \epsilon^c \tau_j^c \!+\! \epsilon^w \!(  \tau_i^{\Delta} \!-\! \tau_j^c\!-\!\tau_j^0 )\big) 
    \!+\!\sum_{ j\in \mathcal{F}_i } \!\epsilon^w \!(  \tau_i^{\Delta}\!-\!\tau_j^0).  \notag
    \end{aligned}
\end{equation}

\vspace{-4.5pt}\subsection{DP-based Platoon and Leader Selection}\vspace{-3pt}
\label{sec:dp final}




\subsubsection{DP recurrence}
In our platoon coordination strategy, we employ DP to decompose the scheduling of \(N\) trucks into a series of sub-problems. 
We define \(V(i), \ i=1, \dots, N\) as the maximum total utility available by scheduling the first~\(i\) trucks, with the final solution given by \(V(N)\)~\cite{bellman1954theory,denardo2012dynamic}.

The key to our method is the DP recurrence relation, which builds the solution for a larger problem by leveraging the optimal results of its sub-problems. In essence, the decision of each truck \(i\) is based on the previously computed optimal value \(V(i-n)\) where $n$ denotes the platoon size, constrained by $1 \! \le \! n \!\le \! \min(i,\bar{n})$ and $\bar{n}$ is the maximum platoon length for safety.
We derive the following recurrence relation, which is the basis of our coordination approach:
\begin{equation}
\label{eq:dp}
V(i) \!=\! 
 \max_{1 \le n \le  \min(i,\bar{n})}  \{
 \max_{{m\in\{\texttt{E}, \texttt{F}\}}} \{ V(i-n) + J(\mathcal{P}_i^{n,m}) \} \},
\end{equation}
where the optimal sub-problem solution for the first $(i-n)$ vehicles has been computed.
The term $J(\mathcal{P}_i^{n,m})$ is the utility of the $i$-th truck in a platoon of size $n$ with leader type $m$, as defined in \eqref{eq:specific_J}.
 
\begin{figure}
    \centering
    \includegraphics[width=0.8\linewidth]{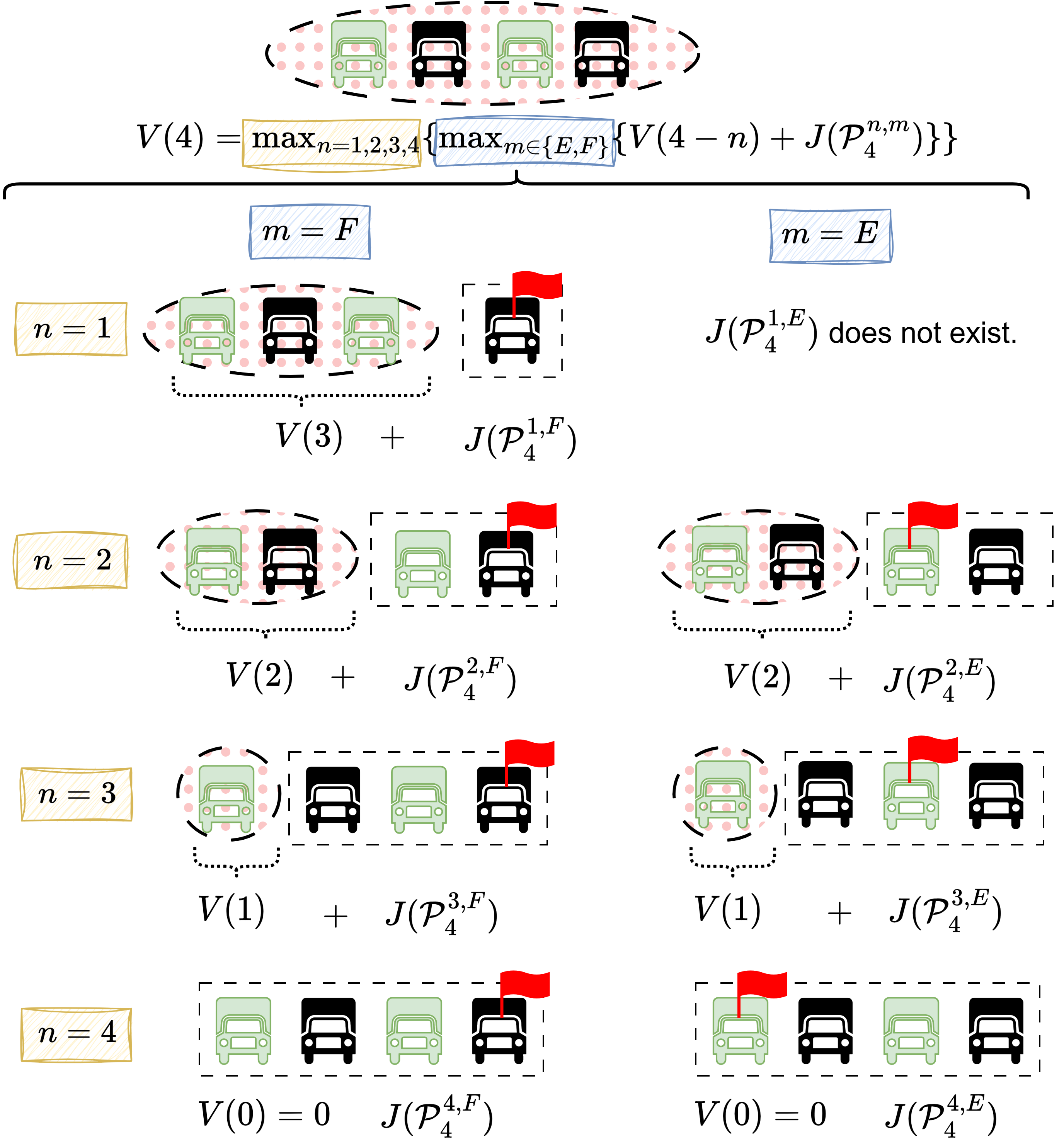}
    \vspace{-7pt}
    \caption{ DP recurrence illustration for $i=4$ to update the optimal platoon utility $V(4)$, where the preceding  optimal values  $V(0),V(1),V(2),V(3)$ have already been obtained. The green trucks represent ETs and the black trucks represent FTs. The leader trucks are marked with red flags. The rectangle illustrates the trucks which will depart in a platoon. The ellipse illustrates the precomputed optimal platoon formation.}
    \label{fig:DP illustration}
    \vspace{-10pt}
\end{figure}

In \eqref{eq:dp}, to determine the optimal platoon formation for the first \( i \) vehicles, where $i\!=\!1, ..., N$, we iterate over all possible platoon sizes and leader choices through two nested maximizations:
(1) External $\max$: We iterate over all possible platoon sizes \( n \), where \( 1 \!\le\! n\!\le\!\min(i,\bar{n}) \).
(2) Internal $\max$: If the platoon consists of both ETs and FTs, we evaluate the utility under both leader options, $J(\mathcal{P}_i^{n,m=\texttt{E}})$ and $J(\mathcal{P}_i^{n,m=\texttt{F}})$, and select the leader type that yields the higher utility to update $V(i)$. If the platoon contains only one type of vehicle, that type is assigned as the leader by default. These two nested maximization steps ensure that all feasible platoon configurations are considered, resulting in optimal platoon formation. Fig. \ref{fig:DP illustration} illustrates how the optimal utility is updated using \eqref{eq:dp}, for the first four trucks in a system with a maximum platoon size \(\bar{n}\!\ge 4\).

Recall that in Section \ref{sec:ordering}, we discretize the solution space under the modest requirement that an ET has enough energy to travel as a follower. Acting as a platoon leader is more demanding with the departure SoC. In accordance with \eqref{eq:tau_c_constraint}, the minimum departure SoC for an leader candidate~\( j \) is $\beta_j v_j d\!+\!s_j^{\text{safe}}$, with $\beta_j\!=\!1$.
If an ET fails to meet this threshold, assigning it as a leader is infeasible.
Thus, we need to verify the SoC requirement when updating $V(i)$, based on the total charging time, \(\tau_{j}^{c}\), obtained from \eqref{eq:charge_for_ev} and the initial SoC, $s_j^0$. This evaluation is termed as \textit{feasibility check}.

\subsubsection{Backtracking}
After sequentially computing the optimal utility $V(i)$ for $i=1, \dots, N$ using \eqref{eq:dp}, the DP-based strategy that yields the utility $V(N)$ is identified as the optimal strategy.
This strategy can be recovered via backtracking, a standard procedure in DP~\cite{denardo2012dynamic}, by tracing the decisions that led to the maximum utility. 
It efficiently reconstructs the optimal platoon schedule and leader assignment.

\begin{algorithm}[t!]
\caption{DP-based Platoon and Leader Selection}
\label{alg:2}
    \textbf{Input:} Loss coefficients \(\epsilon^w, \epsilon^c\), profit coefficients \(\xi^E, \xi^F\);
    for \(i \in \mathcal{N}\): earliest departure time \(\tau_i^{\Delta}\); for \(i \in \mathcal{E}\): departure SoC $s_i^{d,\min}$ and  charging rate \(r_i\).\\
\textbf{Output:} Platoon utility \(V(N)\), platoon formations, leader selections, and departure times.
\begin{algorithmic}[1]
    \STATE Initialize \(V(0)\gets 0\).
    \FOR{\(i=1\) to \(N\)}
        \STATE Set \(V(i) \gets -\infty\)
        \FOR{\(n=1\) to \(\min(i,\bar{n})\)}
            \FOR{\(m \in \{\texttt{E}, \texttt{F}\}\)}
                \STATE Compute candidate utility:
                $
                  U \!\!=\! V(i\!-\!n) \!+\! J(\mathcal{P}_i^{n,m})
                $
                \IF{\(U > V(i)\) and feasibility check is true}
                    \STATE Update \(V(i) \gets U\).
                    \STATE Record the corresponding platoon formation, leader selection, and departure schedule.
                \ENDIF
            \ENDFOR
        \ENDFOR
    \ENDFOR
    \RETURN \(V(N)\) and backtracking the associated strategy.
\end{algorithmic}
\end{algorithm}

\vspace{-4.5pt}\subsection{Algorithm Complexity Analysis}\vspace{-3pt}
The DP-based platoon with leader-selection strategy
is summarized as Algorithm 2. The computational complexity mainly arises from updating the DP value function \(V(i), i\!=\!1,\dots ,N\). To compute \(V(i)\), the algorithm enumerates every possible platoon size \(n\)  with the constraint \(n\!\le\!\bar{n}\) and, for each size, evaluates the two possible leader categories (ET or FT). Thus, updating a single DP state requires at most \(O(\bar{n})\) time. Repeating this computation for all \(N\) vehicles yields a total running time of \(O(N \bar{n})\). When the maximum platoon size is unconstrained, \(\bar{n}\!=\!N\) in the worst case and the algorithm requires \(O(N^{2})\) operations. If \(\bar{n}\) is bounded by a constant, the complexity collapses to \(O(N)\). Thus, although the DP procedure is heuristic and does not guarantee optimality, its time complexity grows at most quadratically with the fleet size $N$ and becomes linear with a constant maximum size \( \bar{n}\), making it computationally efficient for large instances.

\vspace{-3pt}\section{Simulation Study}

\label{sec:Experiments}

\subsection{Parameter Settings}\vspace{-3pt}


We consider a total of $1,000$ trucks in a fleet, consisting of 70\% FTs and 30\% ETs. 
The trucks travel from Stockholm to Linköping at a distance of 200 km, and all of them have a constant speed of 80 km/h.  
We assume that trucks have default departure times from Stockholm, randomly distributed over the full 24-hour period from 00:00 to 23:59. Expressed in minutes, these times range from 1 to 1440. 
The planning horizon is simply set to $T=1440$.
The initial SoCs for ETs are randomly assigned within the range of the safe operation SoC with $s^{\rm safe}_i\!=\!10\%$ to the maximum SoC with $s_i^{\rm max}\!=\!100\%$. We set the maximum platoon size as $\bar{n}\!=\!8$.

Furthermore, we assume that the energy consumption of each follower truck in a platoon is reduced by 18\%, resulting in monetary savings of 0.05€ per km for an ET follower and 0.07€ per km for an FT follower. Given the distance $d\!=\!200$~km, we obtain $\xi^E\!=\!10$€ and $\xi^F\!=\!14$€. 
By the prevailing salary levels of truck drivers in Sweden, the parameter $\epsilon^w$, representing the financial loss associated with waiting, is set at 0.4€ per minute.  
The parameters for ETs are obtained from the most recent published data on Scania-manufactured trucks. Specifically, we consider ETs with a load capacity of 40 tonnes, equipped with a usable battery of 468 kWh, providing a maximum driving range of 350 km. The discharge rate $v_i$ for each truck is set to be 0.286\% SoC per km. The charging power is set at 300 kW, yielding a charging rate $r_i$ of 1.07\% SoC per minute. 
The charge loss $\epsilon^c$ here is set to be 0.2€ per minute. This value is derived from the labor cost of 0.4€ per minute, with a reduction of 0.2€ per minute attributed to a flexibility factor. The flexibility factor reflects the advantage gained by charging, as an increasing SoC enhances the truck’s range and flexibility.

\vspace{-4.5pt}\subsection{Solution Evaluation}\vspace{-3pt}
\begin{figure}
    \centering
    \includegraphics[width=0.8\linewidth]{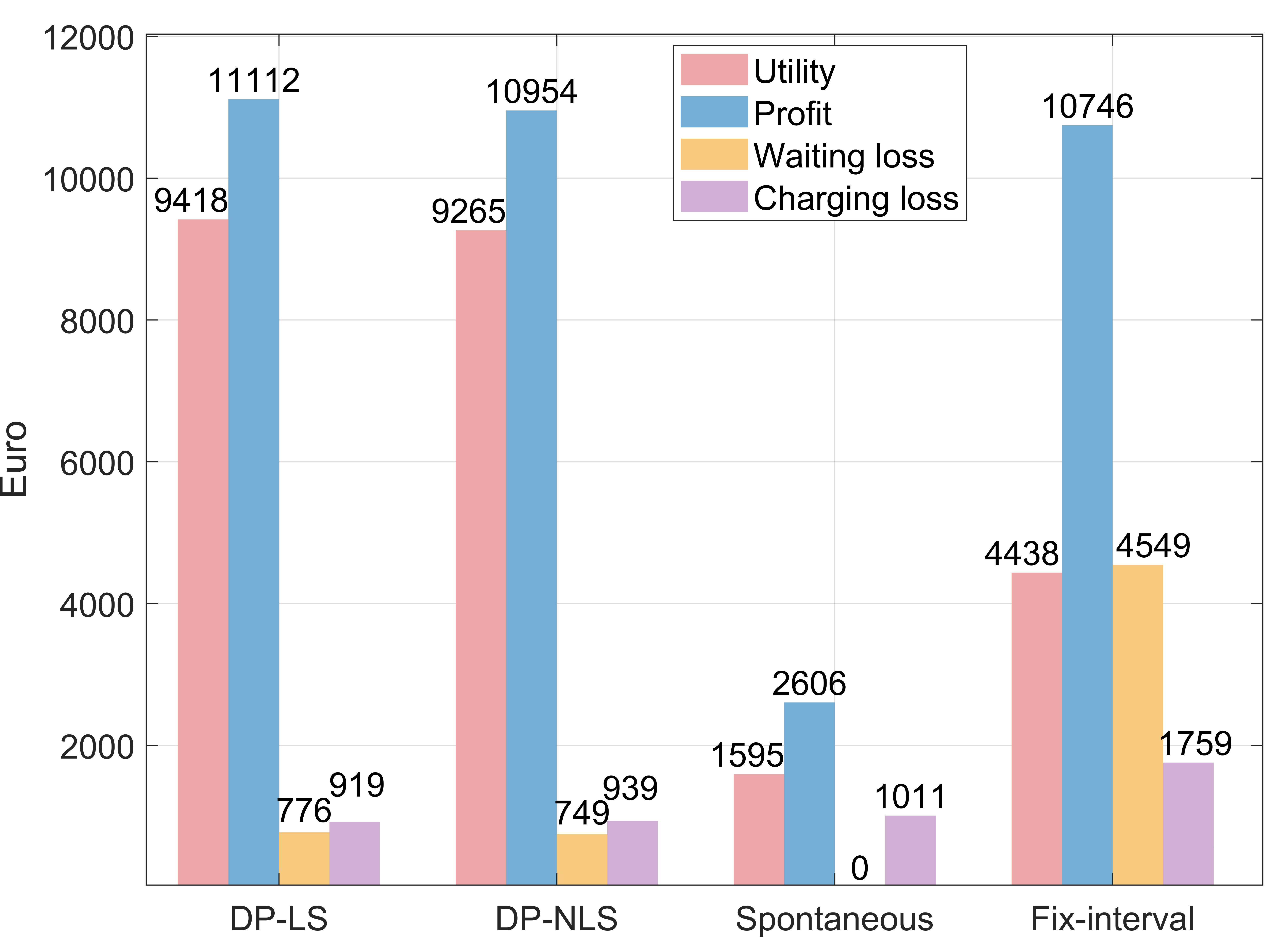}
    \vspace{-10pt}
    \caption{Comparison of the platoon formation results across four methods.}    
    \label{fig:four_cmp}
    \vspace{-10pt}
\end{figure}

We conduct a simulation study to evaluate the platooning performance of four methods: the proposed DP platooning method \textit{with} leader selection (DP-LS),  DP platooning method \textit{without} leader selection (DP-NLS), the spontaneous platooning method, and the fixed-interval platooning method.
In detail, (1) DP-NLS applies DP for departure scheduling but does not consider leader selection, serving as a comparison of DP-LS to assess the effectiveness of leader assignment.
(2) For spontaneous platooning, trucks simply depart at their earliest departure times (default arrival times plus the required ET charging time to safely travel alone) and automatically form platoons with the ones sharing the same departure times. This method does not involve any additional waiting time. 
(3) Fixed‐interval platooning approach partitions the planning horizon into uniform time slots. Trucks whose earliest departure times fall within the same slot are grouped into a platoon and depart simultaneously at the end of that interval. The platoon leader is randomly chosen.

\begin{figure}[t!]
    \centering
    \includegraphics[width=0.8\linewidth]{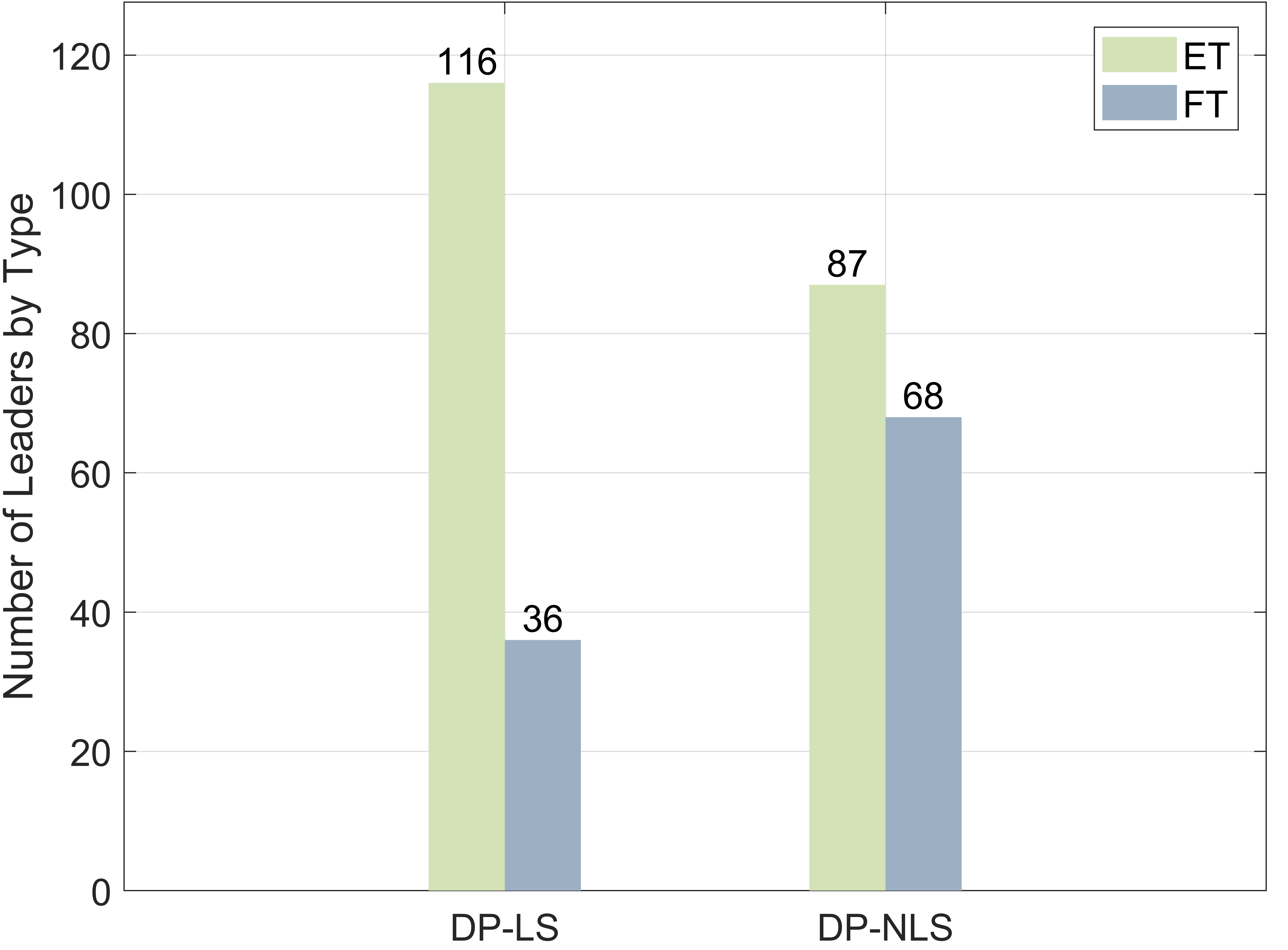}
    \vspace{-10pt}
    \caption{Leader type comparison with DP-LS and DP-NLS methods.}
    \label{fig:ET_FT_distribution}
    \vspace{-3.5pt}
\end{figure}

\begin{figure}[htbp]
    \centering
    \includegraphics[width=0.6\linewidth]{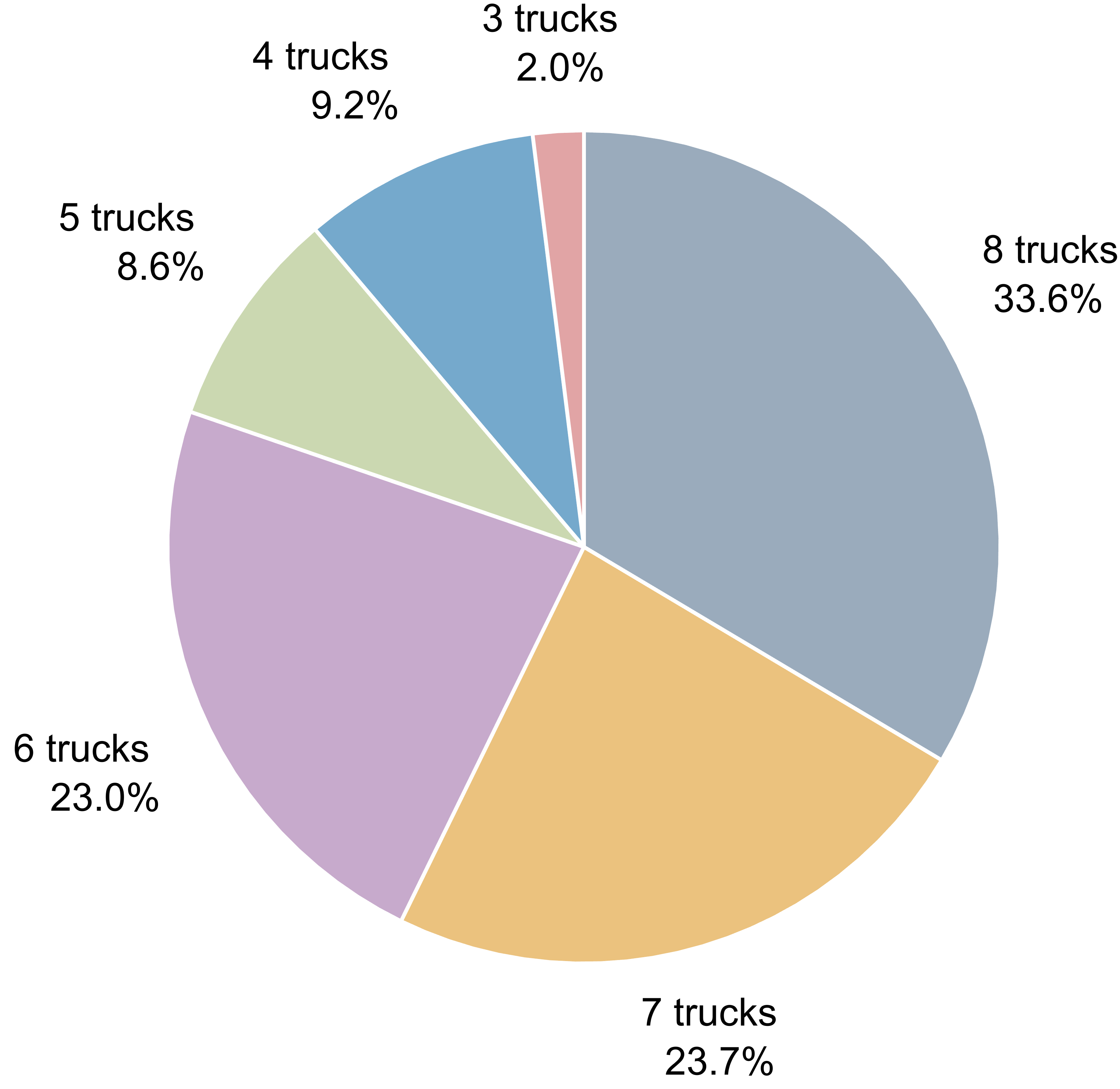}
    \vspace{-10pt}
    \caption{Size distribution of the formed platoons with DP-LS.}
    \label{fig:platoon_size_distribution}
    \vspace{-10pt}
\end{figure}

As shown in Fig.~\ref{fig:four_cmp}, DP-LS method achieves the highest overall utility, demonstrating the effectiveness of leader selection and platoon coordination. DP-NLS, while slightly less efficient, achieves a comparable utility. 
The spontaneous method results in the lowest utility, as it lacks structured coordination. The fixed-interval method also yields a high profit because more trucks can depart in a platoon with a large 30-minute interval, but it also incurs large waiting and charge losses. As a result, its overall utility remains lower despite the potential profits from platooning.

In Fig.~\ref{fig:ET_FT_distribution}, the leader selection result of DP-LS is presented by comparing it with DP-NLS, in which platoon leaders are assigned randomly. The result indicates that with the leader selection, more ETs are chosen as platoon leaders, as an FT follower contributes to a higher profit compared to an ET follower. This is consistent with the utility improvement from DP-NLS to DP-LS in Fig.~\ref{fig:four_cmp}. However, despite the preference for ETs as leaders, 36 platoons are still led by FTs. This outcome is expected given that ETs constitute only 30\% of the $1,000$ trucks in the simulation. The limited availability of ETs limits their roles as platoon leaders.

In this simulation, a total of 152 platoons are formed using DP-LS. The distribution of platoon lengths is illustrated in Fig.~\ref{fig:platoon_size_distribution}. It can be observed that with the proposed method, all trucks are in platoons, and about 80\% of the platoons have a length from 6 to 8, and the left 20\% has a relatively smaller size, varying from 3 to 5. This distribution suggests that DP-LS effectively
groups trucks into well-sized platoons.

In terms of computational efficiency, (1) DP-LS and DP-NLS take around 3.82 seconds. (2) The other two methods take less than 0.01 seconds. While the DP-based methods require higher computational times, they remain acceptable for a large-scale coordination problem. 
Our experiment with $1,000$ trucks confirms this. Moreover, given the achieved improvement in utility, the additional computational cost is justified for practical implementation.

\vspace{-2pt}\section{Conclusion}
\label{sec:Conclusion}
In this paper, we addressed the problem of hub-based truck platoon coordination in a mixed fleet of FTs and ETs. We formulated the problem to jointly optimize ET charging durations, waiting time schedules for both ETs and FTs, and leader assignments to maximize overall platooning utility. To manage the intractable search space while satisfying SoC constraints for ETs with safety guarantees, we introduced earliest departure times to discretize the large solution space and proposed an adjacent platoon formation scheme. We then applied DP to efficiently explore this reduced state space, obtaining suboptimal solutions in polynomial time. Simulation studies with $1,000$ trucks operating between two Swedish hubs validated the proposed approach, demonstrating substantial platooning utility and high computational efficiency. A direction for future research could be deriving analytical bounds on the utility gap between our DP-based approach and the optimum. Another extension is to incorporate nonlinear charging models and potential uncertainties in the trucks' arrival process. 

\bibliographystyle{IEEEtran}

\bibliography{platoon_coordination}

\end{document}